\documentclass{amsart}
\usepackage{amsfonts}
\usepackage{graphicx}
\usepackage{epstopdf}

\newcommand{\beq}{\begin{equation}}
\newcommand{\eeq}{\end{equation}}
\newtheorem{prop}{Proposition}
\newtheorem{thm}{Theorem}
\newtheorem{lemma}{Lemma}
\begin{document}

\title[Universal Casimirs]{Casimir eigenvalues for universal Lie algebra}
       \author{R. L. Mkrtchyan}
       \address{A.I.Alikhanyan National Science Laboratory, 2 Alikhanian Brothers St., Yerevan, 0036, Armenia}
           \email{mrl@web.am}
           
            \author{A. N. Sergeev}
       \address{Department of Mathematics, Saratov State University, Astrakhanskaya 83, Saratov, 410012, Russia}
           \email{sergeevAN@info.sgu.ru}

       \author{A. P. Veselov}
       \address{School of Mathematics,
        Loughborough University, Loughborough,
        Leicestershire, LE11 3TU, UK and
       Moscow State University, Russia}
       \email{A.P.Veselov@lboro.ac.uk}

\maketitle




{\small  {\bf Abstract.} For two different natural definitions of Casimir operators for simple Lie algebras we show that their eigenvalues in the adjoint representation can be expressed polynomially in the universal Vogel's parameters $\alpha, \beta, \gamma$ and give explicit formulae for the generating functions of these eigenvalues. }

\bigskip

\section{Introduction}


Let's introduce a projective plane with coordinates $\alpha, \beta$ and $\gamma$ modulo common multiplier. Quotient space of this projective plane by permutations of these three parameters is called {\it Vogel plane}. This is a moduli space of the tensor category, which is meant to be a model of the {\it universal simple Lie algebra} \cite{Del,V, LM1}. The simple Lie algebras correspond to the values of the parameters given in Table 1, where $t=\alpha+\beta+\gamma.$
\begin{table}[h]  
\caption{Vogel's parameters for simple Lie algebras}     
\begin{tabular}{|c|c|c|c|c|c|}
\hline
Type & Lie algebra  & $\alpha$ & $\beta$ & $\gamma$  & $t=h^\vee$\\   
\hline    
$A_n$ &  $\mathfrak {sl}_{n+1}$     & $-2$ & 2 & $(n+1) $ & $n+1$\\
$B_n$ &   $\mathfrak {so}_{2n+1}$    & $-2$ & 4& $2n-3 $ & $2n-1$\\
$C_n$ & $ \mathfrak {sp}_{2n}$    & $-2$ & 1 & $n+2 $ & $n+1$\\
$D_n$ &   $\mathfrak {so}_{2n}$    & $-2$ & 4 & $2n-4$ & $2n-2$\\
$G_2$ &  $\mathfrak {g}_{2}  $    & $-2$ & $10/3 $& $8/3$ & $4$ \\
$F_4$ & $\mathfrak {f}_{4}  $    & $-2$ & $ 5$& $ 6$ & $9$\\
$E_6$ &  $\mathfrak {e}_{6}  $    & $-2$ & $ 6$& $ 8$ & $12$\\
$E_7$ & $\mathfrak {e}_{7}  $    & $-2$ & $ 8$& $ 12$ & $18$ \\
$E_8$ & $\mathfrak {e}_{8}  $    & $-2$ & $ 12$& $20$ & $30$\\
\hline  
\end{tabular}
\end{table}

The normalisation of parameters follows \cite{LM1} and corresponds to the so-called minimal choice of the bilinear form, when the square of the length of the maximal root is chosen to be 2. In all cases except $C_n$ and $G_2$ this coincides with the choice of Bourbaki, see tables I-VIII at the end of \cite{Bour} (in $C_n$ and $G_2$ cases the Bourbaki's choice is 2 and 3 times the minimal form respectively). Note that the parameter $t=h^\vee$ in this normalisation coincides with the dual Coxeter number.

We are interested in the quantities that can be expressed in terms of the universal parameters in an analytic (say, rational) way.
An example is the dimension of the corresponding simple Lie algebra $\mathfrak {g}$, which can be given by the following expression \cite{Del,V}:
\beq
\label{f3}
dim \, \mathfrak {g} = \frac{(\alpha-2t)(\beta-2t)(\gamma-2t)}{\alpha\beta\gamma}
\eeq

In this note we show that this formula is a consequence of a more general statement:
\beq
\label{new0}
\hat C(z)= - 2 z \frac{d}{dz} \ln \frac{(16t^2-(2t-\alpha)^2 z)(16t^2 -(2t -\beta)^2 z)(16t^2-(2t-\gamma)^2 z)}{(16t^2-\alpha^2 z)(16t^2-\beta^2 z)(16t^2-\gamma^2 z)}
\eeq
for the generating function
$$
\hat C(z)=\sum_{k=0}^{\infty} \hat C_{2k}z^k
$$
of the eigenvalues of certain Casimir operators $\hat C_{2k}$ in the adjoint representation of $\mathfrak {g}$ 
(see Section 3 below). 

We compute also the generating function for the following standard choice of the Casimir operators \cite{Racah, GOR, pp, BR}. Namely, consider the following elements of the centre of the universal enveloping algebra $U \mathfrak g$ 
$$
C_p = g_{\mu_1...\mu_p}X^{\mu_1}...X^{\mu_p}, p=0,1,2,... 
$$
where $X^{\mu}$ are the generators of $\mathfrak g,$ 
$g^{\mu_1...\mu_n}=Tr(\hat{X}^{\mu_1}...\hat{X}^{\mu_n}),$
the trace is taken in the adjoint representation of $\mathfrak g$  and the indices are lowered using the Cartan-Killing form. 
We show that the generating function of the eigenvalues of these Casimir operators in the adjoint representation can be given in terms of the universal parameters by the formula 
\beq
\label{cadinvsym}
C(z)= dim\, \mathfrak g
+  z^2 \frac{16t^3 + 28t^3z + (14 t^3 + t s - 3p)z^2 + (2t^3 + t s - 2p ) z^3}{(2t+\alpha z)(2t+\beta z)(2t+\gamma z)(2+ z)(1+z)}
\eeq
with $dim \, \mathfrak {g}$ is given by (\ref{f3}) and (in Vogel's notations \cite{V}) 
$$s=\alpha\beta + \beta\gamma+\alpha\gamma, \,\, p=\alpha\beta\gamma.$$

All the formulas of this paper work as well for the classical Lie superalgebras $\mathfrak{sl}(m,n)$ and $\mathfrak{osp}(p,q)$ with non-zero $t$ and (super)dimension, which means that $m-n\neq 0, \pm 1$ and $p-q\neq 0,1,2$ respectively. The exceptional cases are not particularly informative since in Vogel's approach $\mathfrak {f}_{4}$ and $\mathfrak {g}_{3}$ are equivalent to $\mathfrak{sl}_3$ and $\mathfrak{sl}_2$ respectively (see \cite{V}) and in the (potentially most interesting) case of $\mathfrak {D}_{2,1,\lambda}$ the parameter $t=\lambda_1+\lambda_2+\lambda_3=0$ (see table below).

 \begin{table}[h]  
\caption{Vogel's parameters for basic classical Lie superalgebras}     
\begin{tabular}{|c|c|c|c|c|}
\hline
 Lie superalgebra  & $\alpha$ & $\beta$ & $\gamma$  & $t$\\   
\hline    
$\mathfrak {sl}_{m,n}$     & $-2$ & 2 & $m-n $ & $m-n$\\
$\mathfrak {osp}_{p,q}$    & $-2$ & 4& $p-q-4 $ & $p-q-2$\\
$\mathfrak {f}_{4}$    & $-2$ & 2& $3$ & $3$\\
$\mathfrak {g}_{3}$    & $-2$ & 2& $2$ & $2$\\
$\mathfrak {D}_{2,1,\lambda}$    & $\lambda_1$ & $\lambda_2$& $\lambda_3$ & $0$\\
\hline  
\end{tabular}
\end{table}

\section{Eigenvalues of universal Casimir operators in adjoint representation}

Following \cite{Racah,GOR} define the universal Casimir operators for any simple Lie algebra $\mathfrak g$ as the following elements of the centre of the corresponding universal enveloping algebra $U \mathfrak g$ 
\beq
\label{nonsymcas}
C_p = g^{\mu_1...\mu_p}X_{\mu_1}\dots X_{\mu_p}, \, p=0,1,2,... 
\eeq
where $X_{\mu}$ are the generators of $\mathfrak g,$ 
\beq
\label{invtens}
g_{\mu_1...\mu_n}=Tr(\hat{X}_{\mu_1}\dots \hat{X}_{\mu_n}),
\eeq
where the trace is taken in the adjoint representation of $\mathfrak g$  and the indices are lifted using the {\it canonical Cartan-Killing form} 
$$
g_{\mu\nu}=Tr(\hat{X}_{\mu}\hat{X}_{\nu})
$$ (see for the details Chapter 9 in \cite{BR}). 

We should mention that some authors (most notably Gelfand  \cite{Gelfand}) use a completely symmetrised version
\beq
\label{symcas}
C^s_p = h^{\mu_1...\mu_p}X_{\mu_1}...X_{\mu_p},
\eeq
where
$$
h_{\mu_1...\mu_n}=\frac{1}{n!}\sum_{\sigma\in S_n}Tr(\hat{X}_{\mu_{\sigma(1)}}\dots \hat{X}_{\mu_{\sigma(n)}})
$$
with the sum taken over all permutations. In this paper we will not consider these symmetric Casimir operators (see though concluding remarks). 

Since adjoint representation is irreducible the Casimir operator $C_p$ acts as a scalar operator with the eigenvalue, which we also denote as $C_p$. We claim that these eigenvalues can be expressed in the terms of the universal parameters $\alpha, \beta, \gamma$ as follows. 

Let as before $$t_2=\alpha^2 + \beta^2+ \gamma^2,\,t_3=\alpha^3+\beta^3+\gamma^3.$$

\begin{thm}
The generating function $C(z)=\sum_{p=0}^{\infty} C_p z^p$ has the form
\begin{align}
\label{cadinv}
&C(z)= \frac{(\alpha-2t)(\beta-2t)(\gamma-2t)}{\alpha\beta\gamma}\\  \nonumber
&+  z^2 \frac{96t^3 + 168t^3z + 6(14 t^3 + t t_2 - t_3)z^2 + ( 13t^3 + 3t t_2 - 4t_3 ) z^3}{6(2t+\alpha z)(2t+\beta z)(2t+\gamma z)(2+ z)(1+z)}.
\end{align}
\end{thm}

In particular, the first few Casimir eigenvalues are
\beq
\label{ucas}
C_0= dim \, \mathfrak {g},\,\,C_1=0,\,\, C_2=1, \,\,C_3=-\frac{1}{4}, \,\,\,  C_4= \frac{3 t  t_2- t_3}{16t^3}.
\eeq

For the proof we use the following expression for $ C_{k}$ found by Okubo \cite{Okubo}:
\beq
\label{cc2}
C_{k}=\sum_{V}\frac{dim V}{dim \, \mathfrak g}\left (\frac{C_2(V)-2C_2(ad)}{2}\right )^k
\eeq
where the sum is taken over all irreducible representations $V$, appearing in the tensor square $ad \otimes ad$ of the adjoint representation $ad$. 

The tensor product $ad \otimes ad$ naturally splits into symmetric and antisymmetric subspaces. The symmetric part in all cases except $\mathfrak {sl}_2$ and $\mathfrak {so}_8$ (which should be treated separately, see below)
according to \cite{Del,V,LM2} decomposes into four  irreducible representations - singlet, and three irreducible representations $Y_2(\alpha), Y_2(\beta), Y_2(\gamma) $ with the dimensions \footnote {For the exceptional Lie algebras we have $3\gamma-2t=0$ (see Table 1), so the component $Y_2(\gamma)$ disappears in this case (cf. \cite{Del}). In the case of  $\mathfrak {sl}_2$ and $\mathfrak {so}_8$ we have $\beta =\gamma$ and the formula (\ref{Y}) can not be used because of the zero at the denominator, see the discussion on page 383 in Landsberg-Manivel  \cite{LM}.}
 \begin{equation}
 \label{Y}
 dim \, Y_2(\alpha)=-\left( \frac{\left( 3\,\alpha - 2\,t \right) \,\left( \beta - 2\,t \right) \,\left( \gamma - 2\,t \right) \,t\,\left( \beta + t \right) \,
      \left( \gamma + t \right) }{\alpha^2\,\left( \alpha - \beta \right) \,\beta\,\left( \alpha - \gamma \right) \,\gamma} \right)
 \end{equation}
(other dimensions given by permutation of parameters). The value of $C_2$ on these representations is zero for singlet, and   $2-\alpha/t, 2-\beta/t, 2-\gamma/t$ respectively (see \cite{V,LM1}), which is actually definition of universal parameters $\alpha, \beta, \gamma$.

Antisymmetric subspace decomposes into the adjoint representation and the remaining part $X_2$  with dimension
$$dim \, X_2= \frac{ dim \, \mathfrak {g} \, (dim \, \mathfrak {g}-3)}{2}$$
and the value of $C_2(X_2)=2$ (see \cite{CM}). Combining all this we have 
\begin{align}\label{lcc}
&C(z)=\frac{1}{dim \, \mathfrak g}\left (\frac{dim \, Y_2(\alpha)}{1+\alpha z/2t}+\frac{dim \,Y_2(\beta)}{1+\beta z/2t}+\frac{dim \, Y_2(\gamma)}{1+\gamma z/2t}+ \frac{1}{1+z}\right)\\ \nonumber
&+\frac{1}{1+z/2}+  \frac{1}{2}(dim\, \mathfrak {g}-3),
\end{align}
which after some calculations leads to the formula  (\ref{cadinv}).
In terms of the symmetric functions
$$s=\alpha\beta + \beta\gamma+\alpha\gamma, \,\, p=\alpha\beta\gamma$$
the generating function has the form (\ref{cadinvsym}).

In the special case of $\mathfrak {sl}_2$ the symmetric part of $ad\otimes ad$ is the sum of trivial representation and 5-dimensional irreducible representation with the eigenvalue $C_2=3.$ 
In the anti-symmetric part we have only adjoint representation. This leads using Okubo's result to the formula
$$C(z)=\frac{5}{3}\frac{1}{1-\frac{1}{2}z}+ \frac{1}{1+\frac{1}{2}z}+\frac{1}{3}\frac{1}{1+z}=3+z^2\frac{3z+4}{(2-z)(2+z)(1+z)}.$$
One can check that this coincides with (\ref{cadinv}) when $\alpha=-2, \beta=\gamma=2.$

A similar calculation for $\mathfrak {so}_8$ gives 
$$C(z)=\frac{1}{28}\left(\frac{300}{1-\frac{1}{6}z}+ \frac{105}{1+\frac{1}{3}z}+\frac{1}{1+z}\right)+ \frac{1}{1+\frac{1}{2}z}+ \frac{25}{2},$$
which coincides with (\ref{cadinv}) for $\alpha=-2, \beta=\gamma=4.$ This completes the proof of Theorem 1.

Direct calculation of the Casimir eigenvalues for unitary and orthogonal groups show the agreement with (\ref{cadinv}). Note that symplectic case is automatically covered by $n \rightarrow -n$ duality \cite{Mkr,Cvitbook,VM}, which is also part of the universality picture.

To show how it works let us consider the eigenvalue of the quartic Casimir $C_4$ for the orthogonal groups $SO(n).$ According to \cite{Cvitbook} (see Table 7.3 on page 72) its value (in a different normalisation) is given by
$$\frac{(n-2)(n^3-9n^2+54n-104)}{8}.$$
The universal parameters of $SO(n)$ can be taken as
$$\alpha=-2, \, \beta=4, \, \gamma=n-4,$$
so $t=n-2.$ So if we assume that the numerator is a symmetric polynomial of $\alpha, \beta, \gamma,$ then we must have an equality
$$n^3-9n^2+54n-104=At^3+Bt t_2 + Ct_3,$$ where 
$$t_2=\alpha^2 + \beta^2+ \gamma^2=n^2-8n+36,\,t_3=\alpha^3+\beta^3+\gamma^3=n^3-12n^2+48n-8.$$
This gives 4 relations on three constants $A,B$ and $C,$ which in general should not be consistent.
In our case however we do have a solution:
$A=0,\, B=3/2, \, C=-1/2$
in agreement with our formula (\ref{ucas}). 

Note that these calculations allow us to predict the eigenvalues of quartic Casimirs for exceptional Lie algebras as well. The comparison with the calculations from \cite{Cvitbook} shows the agreement in all cases except $E_7$, which is probably due to a misprint in \cite{Cvitbook}.

The formula for quartic Casimir operator can be also extracted from the original Vogel's paper \cite{V},
where this operator is denoted as $\pi'$ (see page 17). Theorem 18 of \cite{V}  states that an operator $\sigma$ (which is simply related to $\pi'$) is zero in the adjoint representation ($X_1$ in Vogel's notation). This leads to the expression of $\pi'$ in terms of the universal parameters, which is proportional to our $C_4$. 

It is interesting to note that as $z\rightarrow \infty$ 
we have
\beq
\label{assi}
C(z)\rightarrow (dim \, \mathfrak g-3)/2,
\eeq
which immediately follows from (\ref{lcc}).

\section{Casimir eigenvalues and root systems}

Let  $\mathfrak {h}$ be a Cartan subalgebra of complex simple Lie algebra $\mathfrak {g}$ and $\mathfrak{h}^*$ be its dual space. The root system $R \subset \mathfrak{h}^*$ of $\mathfrak {g}$ is defined as the set of non-zero weights of adjoint representation of $\mathfrak {g}$ with the highest weight  $\lambda_{ad}$ being the {\it maximal root} in $R$ (denoted in \cite{Bour} as $\tilde \alpha$). 

On $\mathfrak {h}$ there is a non-degenerate canonical Cartan-Killing form
\beq
\label{CK}
K(X,Y)= tr\, ad_X ad_Y, \quad X,Y \in \mathfrak {h},
\eeq
where $ad_X: \mathfrak {g} \rightarrow \mathfrak {g}$ is the standard adjoint action defined by $ad_X (Z)= [X, Z].$
In terms of the roots the canonical form can be written as
\beq
\label{CF}
K(X,Y)= \sum_{\alpha \in R} \alpha(X) \alpha (Y)=2 \sum_{\alpha \in R_+} \alpha(X) \alpha (Y)
\eeq
for any choice of positive roots $R_+\subset R.$
The corresponding induced form on $\mathfrak {h}^*$ is denoted in  Bourbaki \cite{Bour} as $\Phi_R$.

The classical result going back to Harish-Chandra says that the algebra of Casimir operators (or the centre of the universal enveloping algebra $U(\mathfrak g)$) is isomorphic  to the algebra of {\it shifted symmetric functions} on $\mathfrak {h}^*$, which are the functions $f$ such that
$f (w\xi-\rho)= f (\xi-\rho), \,\, \xi \in \mathfrak {h}^*$
for every $w \in W$ from the Weyl group $W$, where $\rho \in \mathfrak{h}^*$ is defined in a  usual way as a half-sum of the positive roots
$$\rho = \frac{1}{2} \sum_{\alpha \in R_+} \alpha.$$
In particular, all functions of the form
$f=\phi (\xi +\rho)-\phi(\rho),$ where $\phi \in S^W(\mathfrak h^*)$ is a $W$-invariant function on 
$\mathfrak {h}^*$, have this property.

Consider now the Casimir operator $C^B_{2k}$ corresponding to the function
\beq
\label{CO}
C^B_{2k}(\lambda)=  \sum_{\alpha \in R} [(\lambda + \rho, \alpha)^{2k}-(\rho, \alpha)^{2k}]=2\sum_{\alpha \in R_+} [(\lambda + \rho, \alpha)^{2k}-(\rho, \alpha)^{2k}].
\eeq
Here $\lambda \in  \mathfrak {h}^*$ can be considered as a weight of an irreducible representation $V$ of $\mathfrak g$ with the eigenvalues of $C^B_{2k}$ on $V$ given by (\ref{CO}), where $B=(\, ,\, )$ is some $W$-invariant bilinear form on $\mathfrak {h}^*$ (which is a multiple of the canonical one). 
Note that we can consider here the odd powers as well but the corresponding $C^B_{2k+1}=0.$

The generating function $$C^B(\lambda, z)=\sum_{k=0}^{\infty} C^B_{2k}(\lambda)z^k$$ can be written as
\beq
\label{GF}
C^B(\lambda, z)=- 2z \frac{d}{dz} \ln F^B(\lambda, z), \quad F^B(\lambda, z)= \prod_{\alpha \in R_+} \frac{1-(\lambda + \rho, \alpha)^2z}{1-(\rho, \alpha)^2z}.
\eeq

Now choose as $B$ the {\it minimal invariant bilinear form} on $\mathfrak {h^*}$ normalised by the condition that the maximal root has square 2 (in the induced metric). It is proportional to the canonical form $\Phi_R$:
\beq
\label{dc}
B(\alpha,\beta)=2h^\vee \Phi_R(\alpha,\beta),
\eeq
where $h^\vee$ is the {\it dual Coxeter number} (see \cite{Kac}).

\begin{prop}
For any reduced root system $R$ and the corresponding minimal bilinear form  $B$
\beq
\label{f1}
F^B(\lambda_{ad}, z)=\frac{(4-(2t-\alpha)^2 z)(4-(2t-\beta)^2 z)(4-(2t-\gamma)^2 z)}{(4-\alpha^2 z)(4-\beta^2 z)(4-\gamma^2 z)}
\eeq
with the parameters $\alpha, \beta, \gamma$  given in Table 1.
\end{prop}

To prove this one can use the following key lemma.

\begin{lemma} For any even or odd function $\phi(x)$
\beq 
\label{key}
\prod _{\mu\in R_+}\frac{\phi((\mu,\lambda_{ad}+\rho))}{\phi((\mu,\rho))}=\frac{\phi((\alpha-2t)/2)}{\phi(\alpha/2)}\frac{\phi((\beta-2t)/2)}{\phi(\beta/2)}\frac{\phi((\gamma-2t)/2)}{\phi(\gamma/2)}
\eeq
\end{lemma}

The proof is based on the observation that the roots with $(\mu,\lambda_{ad})\neq 0$ can be
organised in three "strings" with the cancellations between consecutive numerators and denominators within each string (cf. \cite{LM1}).

As a corollary we have a universal formula for the eigenvalues 
\beq
\label{kk}
\hat C_{2k}=\sum_{\alpha \in R} [<\lambda_{ad} + \rho, \alpha>^{2k}-<\rho, \alpha>^{2k}]
\eeq 
of the Casimir operators defined by the canonical form $<\alpha, \beta>=\Phi_R (\alpha, \beta).$
Taking into account the relation (\ref{dc}) and the fact that $t=h^\vee$ for the parameters in Table 1, we have the following

\begin{thm}
The generating function $\hat C(z)=\sum_{i=0}^{\infty} \hat C_{2k}z^k$
of the canonically defined Casimir eigenvalues in the adjoint representation of simple Lie algebras can be expressed in terms of the universal parameters by formula (\ref{new0}). 
\end{thm}

In particular, we have the following universal formulae:
$\hat C_2 = 1$ and 
\beq
\label{ucas1}
\hat C_4 = \frac { 2t^3 + 3t t_2 - t_3}{16 t^3},
\eeq
where $t_2=\alpha^2 + \beta^2+ \gamma^2,\,t_3=\alpha^3+\beta^3+\gamma^3.$
The first formula follows directly from the definition, but the second one is nontrivial (compare it with formula (\ref{ucas}) above).

Note that $$F^B(\lambda_{ad}, z) \rightarrow dim^2\,  \mathfrak g$$
as $z\rightarrow \infty$, which follows from the Weyl expression
for the dimension and reproduces Vogel's formula (\ref{f3}). 
We do not have a similar explanation for (\ref{assi}).

\section{Concluding remarks}

We have considered two different sets of the Casimir operators and found a compact form for the generating functions of their eigenvalues on adjoint representation in terms of the universal parameters. 
 This extends the universal representation of such scalar characteristics of simple Lie (super)algebras as dimensions and the eigenvalues of the quadratic Casimirs  \cite{Del, V, LM1} to the spectra of higher order Casimir operators. 
 
 As we have seen above the very existence of a universal formula of a certain type has a non-trivial consequences and sometimes allows one making calculations completely within classical series to get answers for all  exceptional cases too, which looks remarkable.
 
One can ask about the similar formulas for other Casimir operators, in particular for the minimal degree algebraic generators. The problem is that it seems to be impossible to give a universal definition of such generators for all simple Lie algebras.

An interesting question is to see how all this works for the eigenvalues of the symmetric Casimir operators (\ref{symcas}), which is another universal choice.
In the quartic case these eigenvalues are given in Table 7.3 of \cite{Cvitbook} and  for the orthogonal Lie algebras $so(n)$ are 
$$C_4^s=\frac{(n-2)(n^3-15n^2+138n-296)}{24}.$$
If we make the same assumptions about universal expressions and repeat the calculations as in the non-symmetric case above we can make a prediction (compare with (\ref{ucas}), (\ref{ucas1}))
\beq
\label{pred}
C_4^s= \frac{-4t^3+9t t_2-3 t_3}{48t^3},
\eeq 
 which can be checked for all other simple Lie algebras as well.
Duflo isomorphism \cite{Duflo} could be helpful to compute the universal formulas for higher order symmetric Casimir operators. This possibility is currently under investigation.
 
We would like to mention also that an extension of the universality approach to Chern-Simons theory and affine Kac-Moody algebras is discussed in \cite{MV2}

\section{Acknowledgements}

We thank G. Felder and A.I. Molev for helpful discussions and the referees for their useful comments.

\end{document}